\documentclass[11pt]{article}
\usepackage{cite}
\usepackage{mathrsfs}
\usepackage{amsfonts}
\usepackage{amsmath}
\usepackage{amsfonts,amssymb,color}
\usepackage{dsfont}
\usepackage{curves}
\usepackage{mathrsfs}
\usepackage{pifont}
\usepackage{amssymb}
\usepackage{latexsym,amsmath,amssymb,amsfonts,epsfig,graphicx,cite,psfrag}
\usepackage{eepic,color,colordvi,amscd}

\newtheorem{theorem}{Theorem}[section]

\newtheorem{lemma}{Lemma}[section]

\newtheorem{claim}{Claim}[section]

\newtheorem{observation}{Observation}[section]

\newcommand{\qed}{\hfill\rule{0.5em}{0.809em}}

\def\emptyset{\mbox{{\rm \O}}}

\def\qed{\hfill \rule{4pt}{7pt}}

\def\pf{\noindent {\it Proof. }}

\begin{document}
	
	\title{Coloring of graphs without long odd holes}
	\author{Ran Chen$^{1,}$\footnote{Email: 1918549795@qq.com},  \;Baogang  Xu$^{1,}$\footnote{Email: baogxu@njnu.edu.cn. Supported by 2024YFA1013902},\\\\
		\small $^1$Institute of Mathematics, School of Mathematical Sciences\\
		\small Nanjing Normal University, 1 Wenyuan Road,  Nanjing, 210023,  China\\}
	\date{}
	
	\maketitle
	
	\begin{abstract}
		A {\em hole} is an induced cycle of length at least 4, a $k$-hole is a hole of length $k$, and an  {\em odd hole} is a hole of odd length. Let $\ell\ge 2$ be an integer. Let ${\cal A}_{\ell}$ be the family of graphs of girth at least $2\ell$ and having no odd holes of length at least $2\ell+3$, let ${\cal B}_{\ell}$ be the triangle-free graphs which have no 5-holes and no odd holes of length at least $2\ell+3$, and let ${\cal G}_{\ell}$  be the family of graphs of girth $2\ell+1$ and have no odd hole of length at least $2\ell+5$. Chudnovsky {\em et al.} \cite{CSS2016} proved that every graph in ${\cal A}_{2}$ is 58000-colorable, and every graph in ${\cal B}_{\ell}$  is $(\ell+1)4^{\ell-1}$-colorable. Lan and liu \cite{LL2023} showed that for $\ell\geq3$, every graph in ${\cal G}_{\ell}$ is 4-colorable.  It is not known whether there exists a small constant $c$ such that graphs of ${\cal G}_2$ are $c$-colorable. In this paper, we show that every graph in ${\cal G}_2$ is 1456-colorable, and every graph in ${\cal A}_{3}$ is 4-colorable. We also show that every 7-hole free graph in ${\cal B}_{\ell}$  is $(12\ell+8)$-colorable.
		
		\begin{flushleft}
		\end{flushleft}
		
	\end{abstract}
	
	\newpage
	
	\section{Introduction}
	
	All graphs in this paper are finite and simple. We follow \cite{BM08} for undefined notations and terminologies. Let $G=(V(G), E(G))$ be a graph. Let $v\in V(G)$, and let $X$ and $Y$ be two subsets of $V(G)$. We say that $v$ is {\em complete} to $X$ if $v$ is adjacent to all vertices of $X$, and say that $v$ is {\em anticomplete} to $X$ if $v$ is not adjacent to any vertex of $X$. We say that $X$ is complete (resp. anticomplete) to $Y$ if each vertex of $X$ is complete (resp. anticomplete) to $Y$. We use $G[X]$ to denote the subgraph of $G$ induced by $X$.
	
	We say that a graph $G$ contains a graph $H$ if $H$ is isomorphic to an induced subgraph of $G$. For a set ${\cal H}$ of graphs, we say that $G$ is ${\cal H}$-free if $G$ contains no $H\in {\cal H}$. If ${\cal H}=\{H_1,..., H_t\}$, we simply write $G$ is $(H_1,..., H_t)$-free instead.
	
	Let $u,v\in V(G)$. We simply write $u\sim v$ if $uv\in E(G)$, and write $u\not\sim v$ if $uv\not\in E(G)$. Let $N_G(v)$ be the set of neighbors of $v$, $d_G(v)=|N_G(v)|$. A $uv$-{\em path} is a path connecting $u$ and $v$. The length of a shortest $uv$-path is called the {\em distance} between $u$ and $v$ and denoted by $d_G(u,v)$. For a subset $X$ of $V(G)$ and a vertex $y\in V(G)$, let $d_G(y,X)=\min\{d_G(y,x):~x\in X\}$. Let $i\geq1$ be an integer, and let $N_G^i(u)=\{x\in V(G)~|~d_G(x,u)=i\}$, where $N^1_G(u)=N_G(u)$, which is the set of neighbors of $u$ in $G$. We usually omit the subscript $G$ and simply write $d(v)$, $N(v)$, $d(u,v)$, $d(y,X)$ and $N^i(u)$.
	
	Let $k$ be a positive integer. A $k$-{\em coloring } of $G$ is a function $\phi: V(G)\rightarrow \{1,\cdots,k\}$ such that $\phi(u)\ne \phi(v)$ if $u\sim v$. The {\em chromatic number} $\chi(G)$ of $G$ is the minimum number $k$ for which $G$ has a $k$-coloring. A {\em clique} of $G$ is a set of mutually adjacent vertices in $G$. The {\em clique number} of $G$, denoted by $\omega(G)$, is the maximum size of a clique in $G$. A graph is {\em perfect} if all its induced subgraphs $H$ satisfy $\chi(H)=\omega(H)$. A {\em hole} is an induced cycle of length at least 4, an {\em odd hole} (resp. {\em even hole}) is a hole of odd (resp. even) length, and a $k$-hole is a hole of length $k$. In 2006, Chudnovsky {\em et al} \cite{CRST2006} proved the {\em Strong Perfect Graph Theorem} stating that a graph $G$ is perfect if and only if $G$ does not contain an odd hole or its complement. It is certain that $f(x)=x$ is the binding function of perfect graphs.
	
	In 1975,  Gy\'{a}rf\'{a}s \cite{G75} proposed the following important concept. Let ${\cal G}$ be a family of graphs.  If there exists a function $f$ such that $\chi(G)\leq f(\omega(G))$ for all graphs $G$ in ${\cal G}$, then we say that ${\cal G}$ is $\chi$-{\em bounded}, and call $f$ a {\em binding function} of ${\cal G}$.

	
	For any two positive integers $k$ and $\ell$, Erd\H{o}s \cite{E1959} showed that there exists a graph $G$ with $\chi(G)\geq k$ and without cycles of length less that $\ell$.  This implies that, for ${\cal H}\subseteq$\{triangle, 4-hole, 5-hole, 6-hole, $\cdots$ \}, if the class of ${\cal H}$-free graphs is $\chi$-bounded, then $|{\cal H}|$ must be infinite.
	
	In \cite{G87}, Gy\'{a}rf\'{a}s posed three conjectures: The class of odd hole-free graphs is $\chi$-bounded; for integer $\ell\geq 3$, the class of graphs which does not contain hole of length at least $\ell$ is $\chi$-bounded; for integer $\ell\geq 3$, the class of graphs which does not contain odd hole of length at least $2\ell+1$ is $\chi$-bounded. These three conjectures
	have been confirmed recently \cite{CSS2017,CSSS2020,SS2016}. In 2016, Scott and Seymour \cite{SS2016} showed that $\chi(G)\leq\frac{2^{2^{\omega(G)+1}}}{48(\omega(G)+1)}$ if $G$ is odd hole-free. Then, Chudnovsky, Scott and Seymour \cite{CSS2017} proved that (hole of length at least $\ell$)-free graphs are $\chi$-bounded, and later Chudnovsky, Scott, Seymour and Spirkl \cite{CSSS2020} proved that (odd hole of length at least $\ell$)-free graphs are $\chi$-bounded.
	
	For $\ell\geq2$, let ${\cal G}_\ell$ (resp. ${\cal F}_\ell$) denote the family of graphs which have girth at least $2\ell+1$ and have no odd hole of length at least $2\ell+5$ (resp. $2\ell+3$). In \cite{PZ2014}, Plummer and Zha posed a conjecture claiming that every graph in ${\cal F}_2$ is $3$-colorable. Xu, Yu and Zha \cite{XYZ2017} proved that graphs in ${\cal F}_2$ are 4-colorable. Latter, Wu, Xu and Xu \cite{WXX} extended the result of \cite{XYZ2017} to $\cup _{\ell\ge 2} {\cal F}_\ell$, and proved that all graphs in $\cup _{\ell\ge 2} {\cal F}_\ell$ are 4-colorable, and they posed a new conjecture claiming that graphs in $\cup _{\ell\ge 2} {\cal F}_\ell$ are 3-colorable.
	
	Chudnovsky and Seymour  \cite{CS2023} confirmed the conjecture of Plummer and Zha, and proved that graphs in ${\cal F}_2$ are indeed 3-colorable.   Using the technique of Chudnovsky and Seymour  \cite{CS2023}, Wu, Xu and Xu \cite{WXX2024} proved that all graphs in  ${\cal F}_3$ are 3-colorable. Then, Chen \cite{Chen2024} introduced a new approach and proved that all graphs in
	$\cup _{\ell\ge 5} {\cal F}_\ell$, and
	Wang and Wu \cite{WW2024} applied the approach of Chen \cite{Chen2024} and proved that all graphs ${\cal F}_4$ are 3-colorable.
	
	Notice that ${\cal F}_\ell\subsetneq{\cal G}_\ell$ for all $\ell\ge 2$. Very recently, Lan and Liu \cite{LL2023} showed that all graphs in $\cup_{\ell\ge 3} {\cal G}_\ell$ are 4-colorable, and they also asked that whether every graph in ${\cal G}_2$ is 4-colorable or not. It is even not known wether there is a small constant $c$ such that all graphs in ${\cal G}_2$ are $c$-colorable.
	
	For $\ell\geq2$, let ${\cal A}_{\ell}$ be the family of graphs of girth at least $2\ell$ and having no odd holes of length at least $2\ell+3$ (notice that ${\cal F}_\ell\subsetneq{\cal A}_\ell$ for all $\ell\ge 2$), and let ${\cal B}_{\ell}$ be the triangle-free graphs which have no 5-holes and no odd holes of length at least $2\ell+3$. When $\ell=2$, all graphs in ${\cal A}_{\ell}\cup {\cal G}_{\ell}$ are triangle-free and have no odd holes of length at least 9. The only difference between them is that the graphs in ${\cal A}_{2}$ is 7-hole free but may have 4-holes, and the graphs in ${\cal G}_{2}$ is 4-hole free but may have 7-holes. Furthermore, we see that ${\cal A}_{3}$ consists of exactly those 5-hole free graphs of ${\cal G}_{2}$.
	
	Chudnovsky, Scott and Seymour \cite{CSS2016}  proved that all graphs in ${\cal A}_{2}$ (those
	triangle-free graphs without odd holes of length at least 7) is 58000-colorable. They also showed that all graphs in ${\cal B}_{\ell}$, $\ell\ge 2$, are $(\ell+1)4^{\ell-1}$-colorable.
	
	In this paper, we prove the following two results.
	
	\begin{theorem}\label{main-1}
		Every graph $G$ in ${\cal G}_2$ is $1456$-colorable. If $G$ is further $5$-hole free (i.e., $G\in {\cal A}_{3}$), then $G$ is $4$-colorable.
	\end{theorem}
	
	In \cite{CSS2016}, Chudnovsky, Scott and Seymour pointed out that graphs in ${\cal A}_{2}$ might all be 4-colorable, and they also mentioned that ${\cal A}_{2}$ have non-3-colorable graphs (such as Mycielski-Gr\"ostzsch graph).
	Following the second conclusion of Theorem~\ref{main-1}, we see that graphs in ${\cal A}_{3}$ are all 4-colorable. It might be true that all graphs in $\cup_{\ell\ge 2} {\cal A}_{\ell}$ are 4-colorable.
	
	Our following theorem says that by excluding 7-holes, the chromatic number of graphs of ${\cal B}_{\ell}$ may decrease significantly.
	
	\begin{theorem}\label{main-2}
		Let ${\ell}\geq 2$ be an integer, and let $G$ be a $7$-hole free graph in ${\cal B}_{\ell}$. Then, $\chi(G)\leq 12\ell+8$.
	\end{theorem}

	We will prove the Theorem~\ref{main-1} in Section~\ref{main1}, and prove the Theorem~\ref{main-2} in Section~\ref{main2}.

	\section{Proof of Theorem~\ref{main-1}}\label{main1}
	
	In this section, we will prove the Theorem~\ref{main-1}. Before that, we need to introduce some important notations and prove some useful lemmas.
	
	Let $G$ be a graph. A sequence $(L_0,L_1,\cdots,L_k)$ of disjoint subsets of $V(G)$ is called a {\em levelling} of $G$ if $|L_0|=1$, and for $1\leq i\leq k$, each vertex in $L_i$ has a neighbor in $L_{i-1}$, and has no neighbor in $L_h$ for $0\leq h\leq i-2$; and we say that a levelling $(L_0,L_1,\cdots,L_k)$ is {\em stable} if $L_i$ is stable for $0\leq i\leq k-1$.
	
	For a levelling $(L_0,\cdots,L_k)$ of $G$, and for $0\leq h\leq j\leq k$, we say that a vertex $u\in L_h$ (resp. $v\in L_j$) is an {\em ancestor} (resp. {\em descendant }) of a vertex $v\in L_j$ (resp. $u\in L_h$) if there is a path between $u, v$ of length $j-h$ which has exactly one vertex in each $L_i$ for $h\leq i\leq j$; and if $j=h+1$, we say that $u$ (resp. $v$) is a {\em parent} (resp. {\em child}) of $v$ (resp. $u$). In particular, we say that $v$ is a {\em dependent} of $u$ if $u$ is the unique parent of $v$ in $(L_0,\cdots,L_k)$. For positive integer $i$ and $u, v\in L_i$, an induced $uv$-path with interior in $\cup_{j=0}^{i-1} L_j$ is called a {\em ceiling}-$uv$-{\em path}, and an induced $uv$-path with interior in $\cup_{j=i+1}^{k} L_j$ is called a {\em floor}-$uv$-{\em path}.


	In this section, we say that an odd hole is a {\em long odd hole} if it has length more than seven, and say that a cycle is a {\em short cycle} if it has length less than five. It is certain that each graph in ${\cal G}_2$ has no short cycle and no long odd hole.
	
	\begin{observation}\label{observation}
		Let $G\in {\cal G}_2$, $C=v_1v_2\cdots v_\alpha v_1$ be an odd hole of $G$, and $u\in N(V(C))$. Then the following hold.
		\begin{itemize}
			\item if $\alpha=5$, then $|N_{V(C)}(u)|=1$.
			\item if $\alpha=7$, then either $|N_{V(C)}(u)|=1$ or $N_{V(C)}=\{v_i,v_{i+3}\}$ for some $1\leq i\leq 7$. (the subscript is modulo $7$).
			\item $G$ is bipartite if and only if $G$ is $(5$-hole, $7$-hole$)$-free.
		\end{itemize}
		
	\end{observation}
	
	\begin{lemma}\label{L-1}
		Let $G\in {\cal G}_2$, $(L_0,\cdots,L_k)$ be a stable levelling of $G$, and $H$ be an induced subgraph of $G[L_k]$. Then, $\chi(H[N_H^2(z)])\leq2$ for each $z\in V(H)$.
	\end{lemma}
	\pf We may suppose that $k\ge 2$. Let $z\in L_k$, and let $H'=H[N_H^2(z)]$. We first show that
	$H'$ is 5-hole free. Suppose to its contrary that $H'$ contains a 5-hole, say $C=v_1v_2v_3v_4v_5v_1$. Let $S=\{v_1,\cdots,v_5\}$. Since $S\subseteq N_H^2(z)$, there exists $S'=\{u_1,\cdots, u_5\}\subseteq N_H^1(z)$ such that, for each $i\in \{1,\ldots,5\}$, $u_i$ is  the unique neighbor of $v_i$  in $S$  by Observation~\ref{observation}. Let $z_0$ and $u_4'$ be parents of $z$ and $v_4$ in $L_{k-1}$ respectively. It is certain that $u_4'\not\sim z$ and $z_0$ is anticomplete to $S\cup S'$ as $G$ has no short cycle, and $u_4'\not\sim z_0$ as $L_{k-1}$ is a stable set. By Observation~\ref{observation}, we have that $u_4'$ is anticomplete to $S\setminus\{v_4\}$. If $u'_4\not\sim u_1$, then the induced floor-$z_0u_4'$-path $z_0zu_1v_1v_2v_3v_4u_4'$, together with an even ceiling-$z_0u_4'$-path, produces a long odd hole. Therefore, $u_4'\sim u_1$. By symmetry, we have that $u_4'\sim u_2$, and then a short cycle $zu_1u_4'u_2z$ appears. Therefore, $H'$ is 5-hole-free.

	To prove $\chi(H[N_H^2(z)])\leq2$, it suffices to show further that $H'$ is 7-hole-free (then $H'$ is bipartite). Suppose to its contrary, let $C'=v_1v_2\cdots v_7v_1$ be a 7-hole of $H'$, and let $T=\{v_1,\cdots,v_7\}$. Choose $T_1\subseteq N_H^1(z)$ minimal such that each vertex of $T$ has a neighbor in $T_1$. By Observation~\ref{observation}, we may partition $T_1$ into $T_{11}$ and $T_{12}$ such that $T_{11}=\{u\in T_1~|~|N_{T}(u)|=1\}$ and $T_{12}=\cup_{i=1}^7 \{u\in T_1~|~N_T(u)=\{v_i, v_{i+3}\}$.  Let $z_0$ be a parent of $z$ in $L_{k-1}$. Since $G$ has no short cycle, $z_0$ is anticomplete to $T\cup T_1$, $T_1$ is a stable set, and
	\begin{equation}\label{equ-1}
		\mbox{$N_T(u)\cap N_T(v)=\emptyset$ for any $u, v\in T_1$.}
	\end{equation}
	
	It is certain that $T_{11}\neq\emptyset$. We next prove that
	\begin{equation}\label{equ-2}
		|T_{11}|=1.
	\end{equation}
	
	Suppose to its contrary, let $u$ and $v$ be two vertices of $T_{11}$. Without loss of generality, let $N_T(u)=\{v_1\}$ and $N_T(v)=\{v_i\}$ for some $i$. If $i=3$, then $zuv_1v_7\cdots v_3vz$ is a 9-hole. If $i=6$, then $zuv_1v_2\cdots v_6vz$ is a 9-hole. Therefore, $i\not\in \{3, 6\}$.
	
	Suppose that $i=2$, let $u'_5$ be a parent of $v_5$ in $L_{k-1}$. By Observation~\ref{observation}, $N_T(u'_5)\subseteq\{v_1, v_2, v_5\}$. Since $G$ has no short cycle, we have that $u'_5\not\sim z$ and $u'_5$ is not complete to $\{v_1, v_2\}$. Without loss of generality, we suppose that $u'_5\not\sim v_2$. If $u'_5\not\sim u$, then the floor-$z_0u_5'$-path $z_0zuv_1v_7v_6v_5u'_5$, together with a ceiling-$u'_5z_0$-path, produces a long odd hole. Therefore, $u'_5\sim u$, and hence $u'_5\not\sim v$ to avoid a 4-hole $uzvu'_5u$.  If $u'_5\not\sim v_1$, then the floor-$z_0u_5'$-path $z_0zvv_2v_3v_4v_5u'_5$, together with a ceiling-$u'_5z_0$-path, produces a long odd hole. If $u'_5\sim v_1$ then $uu'_5v_1u$ is a triangle of $G$. Both are contradictions. Therefore, $i\neq 2$, and $i\neq 7$ by symmetry.
	
	Now, we have that $i=4$ or 5. This shows that, for any two vertices $x, y\in T_{11}$, $d_T(x, y)=3$, and so $T_{11}=\{u, v\}$. For $i\in \{5, 6, 7\}$, let $u_i$ be a neighbor of $v_i$ in $T_{12}$. By (\ref{equ-1}), we have that $N_T(u_7)=\{v_7, v_3\}$, which forces $N_T(v)=\{v_6, v_2\}$. But then, $N_T(u_5)\cap N_T(u_6)\neq\emptyset$ or $N_T(u_5)\cap N_T(u_7)\neq\emptyset$, a contradiction to (\ref{equ-1}). Therefore, $|T_{11}|=1$, and (\ref{equ-2}) holds.
	
	\medskip

	By (\ref{equ-1}) and (\ref{equ-2}), we may by symmetry assume that $T_{11}=\{u_1\}$ and $T_{12}=\{u_2,u_3,u_4\}$ such that $N_T(u_1)=\{v_1\}$ and $N_T(u_i)=\{v_i,v_{i+3}\}$ for $2\leq i\leq 4$. Let $u_4'$ be a parent of $v_4$ in $L_{k-1}$. Clearly, $u_4'\not\sim z$ as $G$ has no short cycle. By Observation~\ref{observation}, $N_T(u_4')\subseteq\{v_1,v_4,v_7\}$. If $u_4'$ is anticomplete to $\{v_1,u_1\}$, then the floor-$u_4'z_0$-path $u_4'v_4v_3v_2v_1u_1zz_0$, together with an even ceiling-$u_4'z_0$-path, produces a long odd hole. So, $u_4'$ has exactly one neighbor in $\{v_1,u_1\}$, and $u_4'\not\sim u_2$ to avoid a short cycle $zu_1u_4'u_2z$ or $u_4'v_1v_2u_2u_4'$.
	
	If $u_4'\sim v_1$, then $N_T(u_4')=\{v_1,v_4\}$, and so the floor-$u_4'z_0$-path $u_4'v_1v_7v_6v_5u_2zz_0$, together with an even ceiling-$u_4'z_0$-path, produces a long odd hole. Hence, $u_4'\not\sim v_1$. Now, $u_4'\sim u_1$, and so $u_4'\not\sim v_7$ to avoid a short cycle $u_4'u_1v_1v_7u_4'$.
	
	Let $u_7'$ be a parent of $v_7$ in $L_{k-1}$. Clearly, $u_4'\not=u_7'$. To avoid short cycle $u_7'v_7u_4zu_7'$ or $u_7'v_7v_6u_3u_7'$, we have that $u_7'$ is anticomplete to $\{z,u_3\}$. By Observation~\ref{observation}, $N_T(u_7')\subseteq\{v_3,v_4,v_7\}$. We can deduce that $u_7'\sim v_3$ as otherwise the floor-$u_7'z_0$-path $z_0zu_3v_3v_2v_1v_7u_7'$, together with an even ceiling-$u_7'z_0$-path, produces a long odd hole. So, $N_T(u_7')=\{v_3,v_7\}$.
	
	To avoid a short cycle $u_7'v_3v_2u_2u_7'$, $u_7'\not\sim u_2$. But then, the floor-$u_7'u_4'$-path $u_7'v_7v_1v_2u_2v_5v_4u_4'$, together with an even ceiling-$u_7'u_4'$-path, produces a long odd hole.
	
	Therefore, $H'$ is 7-hole free, and thus $H'$ is bipartite. This proves Lemma~\ref{L-1}. \qed
	
	\medskip
	
	Using a similar argument, we can show that, for any $z\in L_k$, $G[N_{L_k}^3(z)]$ also has an interesting structure.

	\begin{lemma}\label{L-2}
		Let $G\in {\cal G}_2$, $(L_0,\cdots,L_k)$ be a stable levelling of $G$, and $H$ be an induced subgraph of $G[L_k]$. Let $z\in V(H)$, and let $A$ be the set of all vertices of $N_H^3(z)$ whose parents in $L_{k-1}$ are complete to $z$. Then, $\chi(H[A])\leq2$.
	\end{lemma}
	\pf
	Let $H'=H[A]$. We first show that $H'$ is 5-hole free. Suppose to its contrary that $H'$ contains a 5-hole $C=v_1v_2v_3v_4v_5v_1$. Let $v_1v_1'v_1''$ be an induced path such that $v_1'\in N_H^2(z)$ and $v_1''\in N_H^1(z)$. Let $x$ and $y$ be the parents of $v_1''$ and $v_4$ in $L_{k-1}$ respectively. Clearly, $x\not\sim z$ as $G$ is triangle-free. By the definition of $A$, $x$ is anticomplete to $V(C)$ and $y\sim z$. Since $G$ has no short cycle, $x\not\sim v_1'$ and $y$ is anticomplete to $\{v_1',v_1'',v_1,v_2,v_3\}$. Therefore, the union of the floor-$xy$-path $xv_1''v_1'v_1v_2v_3v_4y$ and an even ceiling-$xy$-path is a long odd hole. So, $H'$ is 5-hole free.
	
	Suppose $H'$ contains a 7-hole $C'=v_1v_2\cdots v_7v_1$. Let $S_3=\{v_1,\cdots,v_7\}$. We choose $S_2\subseteq N_H^2(z)$ minimal such that each vertex in $S_3$ has a neighbor in $S_2$, and choose $S_1\subseteq N_H^1(z)$ minimal such that each vertex in $S_2$ has a neighbor in $S_1$. We next show that
	\begin{equation}\label{equ-3}
		\mbox{for any $u\in S_2$, $N_{S_3}(u)=\{v_i,v_{i+3}\}$ for some $1\leq i\leq 7$.}
	\end{equation}
	
	
	By Observation~\ref{observation}, we need only to show that $|N_{S_3}(u)|\geq2$. Suppose not, and we may by symmetry assume that $N_{S_3}(u)=\{v_1\}$. Let $u'$ be a neighbor of $u$ in $S_1$, and let $x$ and $y$ be parents of $u'$ and $v_4$ in $L_{k-1}$ respectively. Clearly, $x\not\sim z$, and so by the definition of $A$, $x$ is anticomplete to $S_3$ and $y\sim z$. Since $G$ has no short cycle, $x\not\sim u$ and $y$ is anticomplete to $\{v_2,v_3\}$. If $y\not\sim v_1$, then the floor-$xy$-path $xu'uv_1v_2v_3v_4y$ can be completed to a long odd hole via an even ceiling-$xy$-path. Hence, $y\sim v_1$. Let $y'$ be a parent of $v_5$ in $L_{k-1}$. With the same argument, we have that $y'$ is complete to $\{z,v_1\}$. But then, there is a short cycle $yv_1y'zy$. This proves (\ref{equ-3}).
	
	Since $G$ has no short cycle, by the definition of $A$, we have that
	\begin{equation}\label{equ-4}
		\mbox{for any $y,y'\in L_{k-1}$, $N_{S_3}(y)\cap N_{S_3}(y')=\emptyset$.}
	\end{equation}
	
	We next prove the following claim.
	\begin{claim}\label{cla-1}
		For any $u\in S_2$ with $N_{S_3}(u)=\{v_i,v_{i+3}\}$ for some $1\leq i\leq7$, there are two vertices $y_{i+3},y_{i+4}\in L_{k-1}$ such that $N_{S_3}(y_{i+3})=\{v_{i+3},v_{i+6}\}$ and $N_{S_3}(y_{i+4})=\{v_{i+4},v_{i+7}\}$.
	\end{claim}
	\pf
	By symmetry, choose $u\in S_2$ with $N_{S_3}(u)=\{v_1,v_{4}\}$. Let $u'$ be a neighbor of $u$ in $S_1$. Let $x$ and $y_4$ be parents of $u'$ and $v_7$ in $L_{k-1}$ respectively. Since $G$ has no short cycle, by definition of $A$, we have that $x$ is anticomplete to $\{z,u\}\cup S_3$, $y_4\sim z$ and $y_4$ is anticomplete to $\{v_5,v_6,u,u'\}$. We can deduce that $y_4\sim v_4$ as otherwise the union of the floor-$xy_4$-path $xu'uv_4v_5v_6v_7y_4$ and an even ceiling-$xy_4$-path is a long odd hole. By Observation~\ref{observation}, $N_{S_3}(y_4)=\{v_4,v_7\}$. Let $y_5$ be a parent of $v_5$. With the same argument, $N_{S_3}(y_5)=\{v_5,v_1\}$. This proves Claim~\ref{cla-1}. \qed
	
	By (\ref{equ-3}), we may by symmetry assume that $u_1\in S_2$ such that $N_{S_3}(u_1)=\{v_1,v_4\}$. By Claim~\ref{cla-1}, there are two vertices $y_4,y_5\in L_{k-1}$ such that $N_{S_3}(y_4)=\{v_4,v_7\}$ and $N_{S_3}(y_5)=\{v_5,v_1\}$. Let $u_2$ be a neighbor of $v_2$ in $S_2$. By (\ref{equ-3}), $N_{S_3}(u_2)\in \{\{v_2,v_5\},\{v_6,v_2\}\}$. If $N_{S_3}(u_2)=\{v_6,v_2\}$, then by Claim~\ref{cla-1}, there is a vertex $y_2\in L_{k-1}$ such that $N_{S_3}(y_2)=\{v_2,v_5\}$; but now $y_2\ne y_5$ and $N_{S_3}(y_2)\cap N_{S_3}(y_5)\ne\emptyset$, which contradicts with (\ref{equ-4}). Therefore, $N_{S_3}(u_2)=\{v_2,v_5\}$.
	
	By Claim~\ref{cla-1}, there is a vertex $y_6\in L_{k-1}$ such that $N_{S_3}(y_6)=\{v_6,v_2\}$. Let $u_3$ be a neighbor of $v_3$ in $S_2$. We have that $N_{S_3}(u_3)\in\{\{v_3,v_6\},\{v_7,v_3\}\}$ by (\ref{equ-3}). If $N_{S_3}(u_3)=\{v_3,v_6\}$, then there is a vertex $y_7\in L_{k-1}$ such that $N_{S_3}(y_7)=\{v_7,v_3\}$; but then $y_7\ne y_4$ and $N_{S_3}(y_4)\cap N_{S_3}(y_7)\ne\emptyset$, which contradicts with (\ref{equ-4}). If $N_{S_3}(u_3)=\{v_7,v_3\}$, then there is a vertex $y_3\in L_{k-1}$ such that $N_{S_3}(y_3)=\{v_3,v_6\}$; but then $y_3\ne y_6$ and $N_{S_3}(y_3)\cap N_{S_3}(y_6)\ne\emptyset$, which contradicts with (\ref{equ-4}).
	
	Therefore, $H'$ is bipartite. This completes the proof of Lemma~\ref{L-2}. \qed

	\begin{lemma}\label{L-3}
		Let $G\in{\cal G}_2$ without 5-hole, $u\in V(G)$, and $L_i=\{x\in V(G)~|~d_G(x,u)=i\}$. Then $G[L_i]$ is bipartite, and hence $G$ is $4$-colorable.
	\end{lemma}
	\pf
	Suppose to its contrary that there exists some integer $k$ such that $G[L_k]$ is not bipartite. We choose the $k$ as small as possible (clearly, $k\geq3$), and so $G[L_k]$ contains a 7-hole $C=v_1v_2\cdots v_7v_1$. For $1\leq i\leq 7$, let $P_i=ux_1^ix_2^i\cdots x_{k-1}^iv_i$ be an induced path from $u$ to $v_i$ of length $k$ such that $x_j^i\in L_j$ for $1\leq j\leq k-1$; since $G$ has no cycle of length at most 5, $x_{k-1}^i$ is anticomplete to $V(C)\setminus\{v_i\}$.
	
	Let $w\in L_h$ be the last common vertex of $P_i$ and $P_{i+2}$, and so $h\leq k-2$. Let $P_i'=wx_{h+1}^i\cdots x_{k-1}^iv_i$ for $1\leq i\leq 7$. There exists a hole $C_1$ such that $V(C_1)\subseteq V(P_i')\cup V(P_{i+2}')\cup \{v_i,v_{i+1},v_{i+2}\}$ and $\{v_i,v_{i+1},v_{i+2}\}\subseteq V(C_1)$. It is certain that $C_1$ is an odd hole as otherwise $G[(V(C_1)\cup V(C))\setminus\{v_{i+1}\}]$ is a long odd hole. Therefore, $C_1$ is a 7-hole, and $x_{k-2}^i\sim x_{k-2}^{i+2}$. We can easily deduce that $G[\{x_{k-2}^1,\cdots,x_{k-2}^7\}]$ contains either a cycle of length at most 5 or a 7-hole, both are contradictions. This proves Lemma~\ref{L-3}. \qed

	\begin{lemma}\label{L-4}
		Let $G\in {\cal G}_2$, $(L_0,\cdots,L_k)$ be a stable levelling of $G$, and $H$ be an induced subgraph of $G[L_k]$. Then, $\chi(H[N_H^3(z)])\leq7$ for each $z\in V(H)$.
	\end{lemma}
	\pf
	Suppose to its contrary that $\chi(H[N_H^3(z)])\geq8$. Let $A$ be the set of all vertices $v\in N_H^3(z)$ such that each parent of $v$ in $L_{k-1}$ is adjacent to $z$. Then, $\chi(H[A])\leq2$ by Lemma~\ref{L-2}. Let $A_1=N_H^3(z)\setminus A$, and now $\chi(H[A_1])\geq 6$.
	
	Let $z_0$ be a parent of $z$ in $L_{k-1}$, and let $A_2=A_1\setminus N(z_0)$. Since $G$ is triangle-free, we have that $\chi(H[N(z_0)])\leq1$, and so $\chi(H[A_2])\geq 5$. By Lemma~\ref{L-3}, we have that $H[A_2]$ contains a 5-hole $C=v_1v_2v_3v_4v_5v_1$. Let $S_3=\{v_1,\cdots,v_5\}$. Choose $S_2\subseteq N_H^2(z)$ minimal such that each vertex in $S_3$ has a neighbor in $S_2$, and choose $S_1\subseteq N_H^1(z)$ minimal such that each vertex in $S_2$ has a neighbor in $S_1$.
	
	By the minimality of $S_2$ and Observation~\ref{observation}, we may assume that $S_2=\{u_1,\cdots,u_5\}$ such that $N_{S_3}(u_i)=\{v_i\}$. Since $G$ has no short cycle, $u_i$ is anticomplete to $\{u_{i-1},u_{i+1}\}$. We have that $H[S_2]$ is not a 5-hole by Lemma~\ref{L-1}, and so by symmetry suppose $u_1\not\sim u_3$. Let $x$ and $x'$ be neighbors of $u_1$ and $u_3$ in $S_1$ respectively. We have that $x=x'$ as otherwise $H[\{z,x,x',u_1,u_3,v_1,v_3,v_4,v_5\}]$ contains either a short cycle or a long odd hole.
	
	Notice that, by the choose of $A_2$, $z_0$ is anticomplete to $S_3$, and each vertex in $S_3$ has a parent in $L_{k-1}$ which is nonadjacent to $z$. Moreover, $z_0$ is anticomplete to $S_1\cup S_2$ as $G$ has no short cycle.
	
	Let $t_3$ be a parent of $v_3$ in $L_{k-1}$ which is nonadjacent to $z$. To avoid a short cycle $t_3v_3u_3xt_3$, $t_3\not\sim x$. By Observation~\ref{observation}, we have that $t_3\sim u_1$ as otherwise the floor-$z_0t_3$-path $z_0zxu_1v_1v_2v_3t_3$ can be completed to a long odd hole via an even ceiling-$z_0t_3$-path. Let $t_1$ be a parent of $v_1$ in $L_{k-1}$ which is nonadjacent to $z$, and by symmetry, we have that $t_1\not\sim x$ and $t_1\sim u_3$.
	
	If $u_5\sim x$, then $u_3\not\sim u_5$, and so with the same argument above, $t_3\sim u_5$; but now a short cycle $t_3u_5xu_1x$ appears. So, $u_5\not\sim x$. Let $y$ be a neighbor of $u_5$ in $S_1$ such that $x\ne y$. If $u_3\not\sim u_5$, then $H[\{z,x,y,u_3,u_5,v_1,v_2,v_3,v_5,\}]$ contains either a short cycle or a long odd hole. So, $u_3\sim u_5$.
	
	To avoid a short cycle $t_3v_3u_3u_5t_3$, $t_3\not\sim u_5$. We have that $t_3\sim y$ as otherwise the union of the floor-$z_0t_3$-path $z_0zyu_5v_5v_4v_3t_3$ and an even ceiling-$z_0t_3$-path is a long odd hole.
	
	Let $t_4$ be a parent of $v_4$ in $L_{k-1}$ which is nonadjacent to $z$. If $t_4$ is anticomplete to $\{x,u_1\}$, by Observation~\ref{observation}, the union of the floor-$z_0t_4$-path $z_0zxu_1v_1v_5v_4t_4$ and an even ceiling-$z_0t_4$-path is a long odd hole. Hence, $t_4$ has exactly one neighbor in $\{x,u_1\}$.
	
	If $t_4\sim u_1$, then $t_4\not\sim y$ to avoid a short cycle $t_4yt_3u_1t_4$. To avoid short cycles $yu_5u_3y$, $zyu_1xz$, $t_4v_4v_3u_3t_4$ and $t_4v_4v_5u_5t_4$, $y$ is anticomplete to $\{u_1,u_3\}$ and $t_4$ is anticomplete to $\{u_3,u_5\}$. But then, the union of the floor-$z_0t_4$-path $z_0zyu_5u_3v_3v_2v_1u_1t_4$ and an even induced ceiling-$z_0t_4$-path is a long odd hole.
	
	Therefore, $t_4\sim x$. To avoid short cycles $t_4yzxt_4$, $xu_3u_5x$, and $t_1u_3u_5yt_1$, $t_4\not\sim y$, $x\not\sim u_5$ and $t_1\not\sim y$. But now, the union of floor-$t_4t_1$-path $t_4xzyu_5v_5v_1t_1$ and an even ceiling-$t_4t_1$-path is a long odd hole. This completes the proof of Lemma~\ref{L-4}. \qed
	
	\medskip
	
	Let $G$ be a graph, $C\subseteq V(G)$ and let $T=t_1t_2\cdots t_k$ be an induced path of $G$, where $k\geq2$. We call a pair $(C, T)$ a {\em lollipop} of $G$ if 1). $V(T)\cap C=\emptyset$; 2). $G[C]$ is connected; and 3). $t_k$ has a neighbor in $C$; and 4). $C$ is anticomplete to $\{t_1,\cdots,t_{k-1}\}$.
	
	For a lollipop $(C,T)$ of $G$, the {\em cleanliness} of $(C, T)$ in $G$ is the maximum $c$ such that $t_1,\cdots, t_c$ all have distance (in $G$) at least three from $C$ (or 0 if $t_1$ has distance
	two from $C$). It follows that the cleanliness is at most $k-2$. We call $t_1$ the {\em end} of the lollipop. If $(C, T)$ and $(C', T')$ are lollipops in $G$, we say the second is a {\em licking} of the first if $C'\subseteq C$, and they have the same end, and $T$ is a subpath of $T'$, and $V(T')\subseteq V(T)\cap C$ (and consequently the cleanliness of $(C', T')$ is at least that of $(C, T)$).
	
	In \cite{CSS2016}, Chudnovsky, Scott and Seymour proved the following two lemmas. There was a error in the data of the second lemma, which we have corrected here.
	
	\begin{lemma}\label{L-5}{\em \cite{CSS2016}}
		Let $h,\kappa\geq0$ be integers. Let $G$ ba a graph such that $\chi(G[N^2(v)])\leq\kappa$ for every vertex $v$; and let $(C,T)$ be a lollipop in $G$, with $\chi(G[C])>h\kappa$. Then there is a licking $(C',T')$ of $(C,T)$, with cleanliness at least $h$ more than the cleanliness of $(C,T)$, and such that $\chi(G[C'])\geq\chi(G[C])-h\kappa$.
	\end{lemma}
	
	\begin{lemma}\label{L-6}{\em \cite{CSS2016}}
		Let $\ell, \kappa$ be positive integers. Let $G$ be a triangle-free graph without odd hole of length at least $2\ell+3$, such that $\chi(G[N^2(v)])\leq\kappa$ for every vertex $v$. Let $(L_0,L_1\cdots,L_k)$ be a levelling of $G$. Then there is a stable levelling $(M_1,M_2,\cdots,M_t)$ of $G$ such that $\chi(G[M_t])\geq\frac{\chi(G[L_k]-(2\ell-1)\kappa}{2}$.
	\end{lemma}
	
	\begin{lemma}\label{L-7}
		Let $G\in{\cal G}_2$, and let $(L_0,\cdots,L_k)$ be a stable levelling of $G$. Then, $\chi(G[L_k])\leq104$.
	\end{lemma}
	\pf
	Suppose to its contrary that let $H=G[L_k]$, and suppose $H$ is connected and $\chi(H)>104$. By Lemma~\ref{L-1}, $\chi(H[N_H^2(z)])\leq2$ for any $z\in V(H)$. Let $u_0\in V(H)$, $L_1'=\{u_0\}$ and $L_i'=\{u\in V(H)~|~d_H(u,u_0)=i\}$ for $i\geq1$. There exists $k'\geq2$ such that $(L_0',\cdots,L_{k'}')$ is a levelling of $H$ and $\chi(H[L_{k'}'])>52$. By Lemma~\ref{L-7}, there exists a stable levelling $(M_0,\cdots,M_t)$ of $H$ such that
	$$\chi(H[M_t])\geq\frac{\chi(H[L_{k'}']-10}{2}>21.$$
	
	We may choose a component $H[M_t']$ of $H[M_t]$ such that $M_t'\subseteq M_t$ and $\chi(H[M_t'])=\chi(H[M_t])$. For $i=t-1,t-2,\cdots,0$, choose $M_i'\subseteq M_i$ minimal such that each vertex in $M_{i+1}'$ has a neighbor in $M_i'$. For $1\leq i\leq t-1$, by the choice of $M_i'$, we have that each vertex in $M_i'$ has a descendant in $M_t'$. Let $H'=H[\bigcup_{i=0}^tM_i']$. By Lemma~\ref{L-1} and~\ref{L-4}, we have that
	\begin{equation}\label{equ-5}
		\mbox{$\chi(H'[N_{H'}^2(z)])\leq2$ and $\chi(H'[N_{H'}^3(z)])\leq7$ for any $z\in V(H')$.}
	\end{equation}
	
	Since $\chi(H'[M_t'])>1$ and $H'$ is triangle-free, we have that $t\geq2$. Let $a_{t-2}\in M_{t-2}'$, and $X$ be set of descendant of $a_{t-2}$ in $M_{t}'$. By (\ref{equ-5}), $\chi(H'[X])\leq2$, and there exists a component $C_1$ of $H'[M_t']-X$ such that
	\begin{equation}\label{equ-6}
		\chi(C_1)\geq \chi(H'[M_t'])-2>19.
	\end{equation}
	
	Since $H'[M_t']$ is connected, there exists a vertex $a_t\in X$ such that $a_t$ has a neighbor in $V(C_1)$. Let $a_{t-1}\in M_{t-1}'$ such that $a_{t-1}$ is complete to $\{a_t,a_{t-2}\}$. Since $\chi(H'[N_{V(C_1)}(a_t)])\leq1$, by (\ref{equ-6}), there exists a component $C_2$ of $C_1-N_{V(C_1)}(a_t)$ such that
	\begin{equation}\label{equ-7}
		\chi(C_2)\geq \chi(C_1)-1>18.
	\end{equation}
	
	There exists a vertex $b_t\in N_{V(C_1)}(a_t)$ such that $b_t$ has a neighbor in $V(C_2)$ because $C_1$ is connected. Let $b_{t-1}$ be a parent of $b_t$ in $M_{t-1}'$ and $b_{t-2}$ be a parent of $b_{t-1}$ in $M_{t-2}'$. Since $b_t\notin X$ and $H'$ has no short cycle, $b_i\ne a_i$ for $i\in\{t-1,t-2\}$, and $b_{t-1}$ is anticomplete to $\{a_{t-1},a_{t-2}\}$.
	
	Since $\chi(H'[N_{C_2}(b_{t-1})])\leq1$, by (\ref{equ-7}), there exists a component $C_3$ of $C_2-N_{V(C_2)}(b_{t-1})$ such that
	\begin{equation}\label{equ-8}
		\chi(C_3)\geq\chi(C_2)-1>17.
	\end{equation}
	
	We now define a vertex $b_t'$ as following: if $b_t$ has a neighbor in $V(C_3)$, then define $b_t'=b_t$; otherwise $N_{V(C_2)}(b_{t-1})\ne\emptyset$ (since $b_t$ has a neighbor in $V(C_2)$), and then let $b_t'\in N_{V(C_2)}(b_{t-1})$ such that $b_t'$ has a neighbor in $V(C_3)$.
	
	By (\ref{equ-5}), $\chi(H'[N_{H'}^2(a_t)])\leq2$ and $\chi(H'[N_{H'}^3(b_{t-1})])\leq7$. By (\ref{equ-8}), there exists a component $C_4$ of $C_3-(N_{H'}^2(a_t)\cup N_{H'}^3(b_{t-1}))$ such that
	\begin{equation}\label{equ-9}
		\chi(C_4)\geq \chi(C_3)-2-7>8=4\cdot2.
	\end{equation}
	
	Since $C_3$ is connected, there exists a lollipop $(V(C_4),T_1)$ of $H'$ with end $b_{t-1}$, where $T_1$ is an induced path joining $b_{t-1}$ and $V(C_4)$, $b_t'\in V(T_1)$, and $V(T_1)\subseteq \{b_{t-1},b_t'\}\cup (V(C_3)\setminus V(C_4))$.  By (\ref{equ-5}), (\ref{equ-9}) and Lemma~\ref{L-5}, there exists a licking $(V(C_5),T_2)$ of $(V(C_4),T_1)$ such that $C_5=H'[V(C_5)]$, the cleanliness of $(V(C_5),T_2)$ is no less than 4, and
	\begin{equation}\label{equ-10}
		\chi(C_5)\geq\chi(C_4)-4\cdot2>0.
	\end{equation}
	
	We may assume that $T_2=p_1p_2p_3p_4\cdots p_m$, where $p_1=b_{t-1}$ and $p_2=b_t'$. Since cleanliness of $(V(C_5),T_2)$ is no less than 4, $m\geq6$ and $d_{H'}(p_i,V(C_5))\geq3$ for $1\leq i\leq4$. Moreover, we have that
	\begin{equation}\label{equ-11}
		\mbox{$d_{H'}(b_{t-1},V(C_5))\geq4$ and $d_{H'}(a_t,V(C_5))\geq3$.}
	\end{equation}
	
	Let $p_{m+1}\in V(C_5)$ such that $p_m\sim p_{m+1}$, and let $c_{t-1}$ be a parent of $p_{m+1}$ in $M_{t-1}'$ and $c_{t-2}$ be a parent of $c_{t-1}$ in $M_{t-2}'$. It is certain that $c_{i}\ne a_{i}$ for $i\in\{t-1,t-2\}$ and $a_{t-2}\not\sim c_{t-1}$ as $p_{m+1}\notin X$. By (\ref{equ-11}), we have that $c_{i}\ne b_{i}$ for $i\in\{t-1,t-2\}$, $c_{t-2}\not\sim b_{t-1}$, $b_{t-2}\not\sim c_{t-1}$ and $a_t\not\sim c_{t-1}$. Notice that $(M_0',\cdots,M_t')$ is a stable levelling of $H'$, and thus $\{a_i,b_i,c_i\}$ is a stable set for $i\in\{t-1,t-2\}$.
	
	Choose the minimum $1\leq h\leq m+1$ such that $c_{t-1}\sim p_h$. Since $d_{H'}(c_j,V(C_5))\geq3$ for $1\leq j\leq4$, $h\geq5$. Let $P=c_{t-1}p_hp_{h-1}\cdots p_2$ (Notice that $p_2=b_t'$). Clearly, $P$ is an induced path of length at least 4. The union of $P$, $p_2p_1b_{t-2}$, $c_{t-1}c_{t-2}$ and an even ceiling-$c_{t-2}b_{t-2}$-path (in $\cup_{i=0}^{t}M_i'$) is a hole of length at least 9, and so it is an even hole. This implies that $P$ is an induced path of odd length at least 5.
	
	If $b_t'\ne b_t$, then the union of $P$ and $p_2p_1b_ta_ta_{t-1}c_{t-2}c_{t-1}$ is a long odd hole whenever $c_{t-2}\sim a_{t-1}$, and the union of $P$, $p_2p_1b_ta_ta_{t-1}a_{t-2}$, $c_{t-1}c_{t-2}$ and an even ceiling-$a_{t-2}c_{t-2}$-path (in $\cup_{i=0}^{t}M_i'$) is a long odd hole whenever $c_{t-2}\not\sim a_{t-1}$. Therefore, $b_t'= b_t$.
	
	If $c_{t-2}\sim a_{t-1}$, then the union of $P$ and $p_2a_ta_{t-1}c_{t-2}c_{t-1}$ is a long odd hole. If $c_{t-2}\not\sim a_{t-1}$, then the union of $P$, $p_2a_ta_{t-1}a_{t-2}$, $c_{t-1}c_{t-2}$ and an even ceiling-$a_{t-2}c_{t-2}$-path (in $\cup_{i=0}^{t}M_i'$) is a long odd hole.
	
	This completes the proof of Lemma~\ref{L-7}. \qed
	
	\medskip
	
	Now, we can prove Theorem~\ref{main-1}.
	
	\noindent\textbf{{\em Proof of Theorem}~\ref{main-1}}: The second conclusion of Theorem~\ref{main-1} follows from Lemma~\ref{L-3}. Suppose the first conclusion of Theorem~\ref{main-1} does not hold, and we may assume that $G\in {\cal G}_2$ such that $G$ is connected and $\chi(G)>1456$. Therefore, there exists a levelling $(L_0,\cdots,L_k)$ of $G$ such that $\chi(G[L_k])>728$. By Lemma~\ref{L-7}, we have that $\chi(G[N^2(v)])\leq104$ for any $v\in V(G)$. And so, by Lemma~\ref{L-6}, there exists a stable levelling $(M_0,\cdots,M_t)$ of $G$ such that
	$$
	\chi(G[M_t])\geq\frac{\chi(G[L_k])-5\cdot104}{2}>104,
	$$
	which contradicts with Lemma~\ref{L-7}. This completes the proof of Theorem~\ref{main-1}. \qed

	\section{Proof of Theorem~\ref{main-2}}\label{main2}
	
	In this section, we will prove the Theorem~\ref{main-2}. Recall that for a levelling $(L_0,\cdots,L_k)$ of a graph $G$, and two vertices $u,v$ in the levelling, we say that $v$ is a {\em dependent} of $u$ if $u$ is the unique parent of $v$ in $(L_0,\cdots,L_k)$. It is easily to see that for any vertex $v$ in a levelling $(L_0,\cdots,L_k)$ of $G$, if $v$ has no dependent in $(L_0,\cdots,L_k)$, then the sequence obtained from $(L_0,\cdots,L_k)$ by deleting the vertex $v$ is also a levelling of $G$.
	
	We call a levelling $(L_0,\cdots,L_k)$ of $G$ a {\em week stable levelling} if $L_i$ is a stable set for $1\leq i\leq k-2$. Throughout this section, we say that an odd hole is a {\em long odd hole} if it has length at least $2\ell+3$. Before proving Theorem~\ref{main-2}, we need first prove the following useful lemma.
	
	\begin{lemma}\label{L-8}
		Let $\ell\geq2$ be an integer. Let $G$ be a triangle-free graph with no 5-hole, and no odd hole of length at least $2\ell+3$. Let $(L_0,L_1\cdots,L_k)$ be a levelling of $G$. Then there is a week stable levelling $(M_0,M_1,\cdots,M_t)$ of $G$ such that $\chi(G[M_t])\geq\frac{\chi(G[L_k])}{2}-\ell+1$.
	\end{lemma}
	\pf
	We prove the Lemma by induction on $|V(G)|$. If $G$ is bipartite, then the result holds easily, and so we may assume that $\ell\geq3$. If $\chi(G[L_k])\leq 2\ell-2$, then $\frac{\chi(G[L_k])}{2}-\ell+1\leq0$, and so the stable levelling $(L_0,L_1)$ satisfies the lemma. Therefore, we may assume that $\chi(G[L_k])>2\ell-2$. Moreover, by inductive hypothesis, we may assume that $G[L_k]$ is connected, $V(G)=\bigcup_{i=0}^kL_i$, and for $0\leq i\leq k-1$ and every $v\in L_i$, $v$ has a dependent in $L_{i+1}$.
	
	Since $G$ is a triangle-free graph without 5-hole, we can deduce that $\chi(G[N^2(z)])\leq1$ for any $z\in V(G)$.
	
	Let $\ell_0\in L_0$, and $\ell_i\in L_i$ for $1\leq i\leq k$, such that $\ell_i$ is a dependent of $\ell_{i-1}$. Let $L=\{\ell_0,\cdots,\ell_k\}$. Notice that $(L_k,\ell_{k-2}\ell_{k-1})$ is a lollipop of $G$. Since $\chi(G[L_k])>2\ell-2$, by Lemma~\ref{L-5}, there exists a licking $(C_1,T_1)$ of $(L_k,\ell_{k-2}\ell_{k-1})$ such that its cleanliness is no less than $2\ell-2$ and
	\begin{equation}\label{equ-12}
		\chi(G[C_1])\geq\chi(G[L_k])-(2\ell-2)>0.
	\end{equation}
	
	Suppose that $T_1=t_1t_2\cdots t_{2\ell-2}\cdots t_m$, where $m\geq 2\ell$, $t_1=\ell_{k-2}$ and $t_2=\ell_{k-1}$. We have that
	\begin{equation}\label{equ-13}
		\mbox{$d_G(t_i,C_1)\geq3$ for $1\leq i\leq2\ell-2$.}
	\end{equation}
	
	Since $G$ is triangle-free and $\ell_{i+1}$ is a dependent of $\ell_i$, we have that each vertex $v\in N(L)\cap L_i$ is adjacent to exactly one of $\ell_i$ and $\ell_{i-1}$, and has no other neighbors in $L$. We say the {\em type} of a vertex $v\in N(L)\cap L_i$ is $\alpha$, where $\alpha=1$ or 2, depending whether $v$ is adjacent to $\ell_{i-1}$ and not to $\ell_i$, or adjacent to $\ell_i$ and not to $\ell_{i-1}$.
	
	Let $\alpha\in\{1,2\}$, we define $V_\alpha$ be the minimal subset of $V(G)\setminus L$ such that: 1). every vertex in $N(L)$ of type $\alpha$ belongs to $V_\alpha$; and 2). for every vertex $v\in V(G)\setminus(L\cup N(L))$, if some parent of $v$ belongs to $V_\alpha$, then $v\in V_\alpha$.
	
	It is certain that $V(G)=\bigcup_{i=0}^kL_i=L\cup V_1\cup V_2$. Since $\ell_k\sim \ell_{k-1}=t_2$, by (\ref{equ-13}), $\ell_k\notin C_1$. Hence, there exists $\alpha\in \{1,2\}$ such that $\chi(G[V_\alpha\cap C_1])\geq\frac{\chi(G[C_1])}{2}$. Therefore, there is a component $G[C_2]$ of $V_\alpha\cap C_1$, where $C_2\subseteq C_1$, such that
	\begin{equation}\label{equ-14}
		\chi(G[C_2])=\chi(G[V_\alpha\cap C_1])\geq \frac{\chi(G[L_k])}{2}-\ell+1
	\end{equation}
	by (\ref{equ-12}). Since $G[C_1]$ is connected, there is a licking $(C_2,T_2)$ of $(C_1,T_1)$. We may suppose $T_2=t_1t_2\cdots t_{2\ell-2}\cdots t_n$, where $n\geq m$.
	
	Let $J_k=C_2$, and for $i=k-1,k-2,\cdots,1$, choose $J_i\subseteq V_\alpha\cap L_i$ minimal such that each vertex in $J_{i+1}\setminus N(L)$ has a parent in $J_i$. By (\ref{equ-13}), $J_{k-1}$ is anticomplete to $\{t_1,t_2,\cdots,t_{2\ell-2}\}$.
	
	By the choice of $J_i$, we have that for each vertex $v\in J_i$, there exists a vertex $u\in J_{i+1}\setminus N(L)$ such that $v$ is the only parent in $J_i$. Furthermore, for each vertex $p_i\in J_i$, there is an induced path $p_ip_{i+1}\cdots p_k$ such that: 1). $p_j\in J_j$ for $i\leq j\leq k$; 2). $p_j\notin N(L)$ for $i<j\leq k$; and 3). $p_{j-1}$ is the only parent of $p_j$ in $J_{j-1}$ for $i<j\leq k$.  We call such a path a {\em pillar} for $p_i$.
	
	\begin{claim}\label{cla-2}
		For $1\leq i\leq k-2$, if $v\in J_i$ and $v\not\sim \ell_i$, then there exists an induced floor-$v\ell_i$-path $P_v$, of length at least $2\ell-4+2(k-i)$, and for any $u\in J_i\setminus\{v\}$, $u$ has no neighbor in the interior of $P_v$.
	\end{claim}
	\pf
	Let $P=vp_{i+1}p_{i+2}\cdots p_{k}$ be a pillar for $v$. Since $p_k\in J_k$, there is an induced path between $t_2$ and $p_k$ with interior in $L_k$ containing all of $t_2,t_3,\cdots,t_n$. Since $p_{k-1}\sim p_k$, and $p_{k-1}$ is anticomplete to $\{t_2,t_3,\cdots,t_{2\ell-2}\}$, there is an induced path between $t_2$ and $p_{k-1}$ with interior in $L_k$ containing all of $t_2,t_3,\cdots,t_{2\ell-2}$ and at least one more vertex of $L_k$, and thus with length at least $2\ell-2$. Its union with the induced path $vp_{i+1}\cdots p_{k-1}$ and the induced path $\ell_i\ell_{i+1}\cdots \ell_{k-1}$ is an induced path between $v$ and $\ell_{i}$, of length at least $2\ell-4+2(k-i)$. We call such a path $P_v$.
	
	Let $u\in J_i\setminus\{v\}$. Since $i\leq k-2$, $u$ is anticomplete to $L_k$. By the choice of the pillar $P$ and $L$, $u$ is anticomplete to $\{\ell_{i+1}, p_{i+1}\}$. Therefore, $u$ is anticomplete to the interior of $P_v$ as $(L_0,\cdots,L_k)$ is a levevlling of $G$. This proves Claim~\ref{cla-2}. \qed
	
	\begin{claim}\label{cla-3}
		For $1\leq i\leq k-2$, $v\in J_i$, there is an induced path $R_v=vv_{i-1}\cdots v_h\ell_{h-1}\ell_h\ell_{h+1}\cdots\ell_{i}$ if $\alpha=1$, and $R_v=vv_{i-1}\cdots v_h\ell_h\ell_{h+1}\cdots\ell_{i}$ otherwise, where $v_j\in J_j$ for $h\leq j\leq i-1$.
	\end{claim}
	\pf
	For $1\leq i\leq k-2$ and every $v\in J_i$, either $v\in N(L)$ or it has a parent in $J_{i-1}$. Hence, there exists an induced path $v_{i}v_{i-1}\cdots v_{h}$, where $v=v_i$, for some $h\leq i$ such that $v_h\in J_h\cap N(L)$ and $v_j\in J_j\setminus N(L)$ for $h+1\leq j\leq i$. Since $v_h$ has a neighbor in $L$, we have that if $\alpha$=1, then $v_iv_{i-1}\cdots v_h\ell_{h-1}\ell_h\ell_{h+1}\cdots\ell_{i}$ is an induced path; if $\alpha$=2, then $v_iv_{i-1}\cdots v_h\ell_h\ell_{h+1}\cdots\ell_{i}$ is an induced path. We choose some such path and call it $R_v$. It is certain that the path $R_v$ has even length if $\alpha=1$, and odd length otherwise. This proves Claim~\ref{cla-3}. \qed
	
	\begin{claim}\label{cla-4}
		For $1\leq i\leq k-2$, $J_{i}$ is a stable set.
	\end{claim}
	\pf
	Suppose to its contrary that let $u,v\in J_i$ such that $u\sim v$. Since $G$ is triangle-free and $u,v$ have the same type, we have that not both $u,v\in N(L)$, and so we may assume that $v\notin N(L)$. By Claim~\ref{cla-2}, we choose an induced floor-$v\ell_i$-path $P_v$, which has length of at least $2\ell-4+2(k-i)$, and $u$ has no neighbor in interior of $P_v$. By Claim~\ref{cla-3},  if $\alpha=1$, then $R_v=vv_{i-1}\cdots v_h\ell_{h-1}\ell_h\ell_{h+1}\cdots\ell_{i}$, and otherwise, $R_v=vv_{i-1}\cdots v_h\ell_h\ell_{h+1}\cdots\ell_{i}$, where $h\leq i-1$ and $v_j\in J_j$ for $h\leq j\leq i-1$.
	
	It is certain that if $\alpha=1$, then $R_v$ has even length at least 4; and if $\alpha=2$, then $R_v$ has odd length at least 3. Therefore the union of $P_v$ and $R_v$ is a hole of length at least $2\ell+3$, and so it is an even hole. We have that
	\begin{equation}\label{equ-15}
		\mbox{ if $\alpha=1$, $P_v$ has even length; and otherwise $P_v$ has odd length.}
	\end{equation}
	
	By Claim~\ref{cla-3}, if $\alpha=1$, then $R_u=uu_{i-1}\cdots u_{h'}\ell_{h'-1}\ell_{h'}\ell_{h'+1}\cdots\ell_{i}$, and otherwise, $R_u=uu_{i-1}\cdots u_{h'}\ell_{h'}\ell_{h'+1}\cdots\ell_{i}$, where $u_j\in J_j$ for $h'\leq j\leq i-1$. Therefore, by (\ref{equ-15}), the union of $P_v$, $vu$ and $R_u$ is a long odd hole. This proves Claim~\ref{cla-4}. \qed
	
	If $\alpha=1$, then let $M_0=\{\ell_0\}$, $M_i=\{\ell_i\}\cup J_i$ for $1\leq i\leq k$. If $\alpha=2$, then let $M_0=\{\ell_1\}$, $M_i=\{\ell_{i+1}\}\cup J_i$ for $1\leq i\leq k-1$, and $M_k=J_k$. By (\ref{equ-14}), we have that $\chi(G[M_k])\geq \chi(G[J_k])=\chi(G[C_2])\geq \frac{\chi(G[L_k])}{2}-\ell+1$. By Claim~\ref{cla-4}, it is easy to verify that $(M_0,\cdots,M_k)$ is a week stable levelling of $G$. This completes the proof of Lemma~\ref{L-8}. \qed
	
	\medskip
	
	Now, we can prove the Theorem~\ref{main-2}.
	
	\noindent\textbf{{\em Proof of Theorem}~\ref{main-2}}: Suppose the Theorem~\ref{main-2} does not hold. We may assume that $G$ is a connected triangle-free graph with no odd hole of length at most $7$, and no odd hole of length at least $2{\ell}+3$ satisfies that $\chi(G)> 12\ell+8$. There exists a levelling $(L_0,\cdots,L_k)$ of $G$ such that $\chi(G[L_k])\geq \frac{\chi(G)}{2}>6\ell+4$. By Lemma~\ref{L-8}, there exists a week stable levelling $(M_0,M_1,\cdots,M_t)$ of $G$ such that \begin{equation}\label{equ-16}
		\mbox{$\chi(G[M_t])\geq\frac{\chi(G[L_k])}{2}-\ell+1>2\ell+3$.}
	\end{equation}
	
	We may choose $(M_0,\cdots,M_t)$ such that $G[M_t]$ is connected and for each $v\in \bigcup_{i=0}^{t-1}M_i$ has a descendant in $M_t$ (with the same operation in the proof of Lemma~\ref{L-7}). Let $G'=G[\bigcup_{i=0}^tM_i]$. Since $G'$ is a triangle-free graph with no 5-hole and 7-hole, we have that
	\begin{equation}\label{equ-17}
		\mbox{$\chi(G'[N_{G'}^i(v)])\leq 1$ for each $v\in V(G')$ and $1\leq i\leq3$.}
	\end{equation}
	
	By (\ref{equ-16}), $\chi(G'[M_t])>1$, and so $t\geq2$. Let $a_{t-2}\in M_{t-2}$, and $X$ be set of descendant of $a_{t-2}$ in $M_{t}$. By (\ref{equ-17}), $\chi(G'[X])\leq1$, and there exists a component $C_1$ of $G'[M_t]-X$ such that
	\begin{equation}\label{equ-18}
		\chi(C_1)\geq \chi(G'[M_t])-1>2\ell+2.
	\end{equation}
	
	There exists a vertex $a_t\in X$ which has a neighbor in $V(C_1)$ because $G'[M_t]$ is connected. Let $a_{t-1}\in M_{t-1}$ such that $a_{t-1}$ is complete to $\{a_t,a_{t-2}\}$. Since $\chi(G'[N_{V(C_1)}(a_t)])\leq1$, by (\ref{equ-18}), there exists a component $C_2$ of $C_1-N_{V(C_1)}(a_t)$ such that
	\begin{equation}\label{equ-19}
		\chi(C_2)\geq \chi(C_1)-1>2\ell+1.
	\end{equation}
	
	There exists a vertex $b_t\in N_{V(C_1)}(a_t)$ which has a neighbor in $V(C_2)$ because $C_1$ is connected. Let $b_{t-1}$ be a parent of $b_t$ in $M_{t-1}$ and $b_{t-2}$ be a parent of $b_{t-1}$ in $M_{t-2}$. Since $b_t\notin X$, $b_i\ne a_i$ for $i\in\{t-1,t-2\}$ and $a_{t-2}\not\sim b_{t-1}$. Since $G'$ is a triangle-free graph with no 5-hole and $(M_0,\cdots,M_t)$ is a week stable levelling of $G'$, $b_{t-2}$ is anticomplete to $\{a_{t-1},a_{t-2}\}$, $a_{t-1}\not\sim b_t$, and $b_{t-1}\not\sim a_t$.
	
	Since $\chi(G'[N_{C_2}(b_{t-1})])\leq1$, by (\ref{equ-19}), there exists a component $C_3$ of $C_2-N_{V(C_2)}(b_{t-1})$ such that
	\begin{equation}\label{equ-20}
		\chi(C_3)\geq\chi(C_2)-1>2\ell.
	\end{equation}
	
	With the same operation in the proof of Lemma~\ref{L-7}, we define a vertex $b_t'$ as following: if $b_t$ has a neighbor in $V(C_3)$, then define $b_t'=b_t$; otherwise $N_{V(C_2)}(b_{t-1})\ne\emptyset$, and then let $b_t'\in N_{V(C_2)}(b_{t-1})$ such that $b_t'$ has a neighbor in $V(C_3)$.
	
	By (\ref{equ-17}), $\chi(G'[N_{G'}^2(a_t)])\leq1$ and $\chi(G'[N_{G'}^3(b_{t-1})])\leq1$. By (\ref{equ-20}), there exists a component $C_4$ of $C_3-(N_{G'}^2(a_t)\cup N_{G'}^3(b_{t-1}))$ such that
	\begin{equation}\label{equ-21}
		\chi(C_4)\geq \chi(C_3)-2>2\ell-2.
	\end{equation}
	
	There exists a lollipop $(V(C_4),T_1)$ of $G'$ with end $b_{t-1}$ such that $T_1$ is an induced path joining $b_{t-1}$ and $V(C_4)$, $b_t'\in V(T_1)$, and $V(T_1)\subseteq \{b_{t-1},b_t'\}\cup (V(C_3)\setminus V(C_4))$ because $C_3$ is connected. By (\ref{equ-17}), (\ref{equ-21}) and Lemma~\ref{L-5}, there exists a licking $(V(C_5),T_2)$ of $(V(C_4),T_1)$ with cleanliness at least $2\ell-2$, such that $C_5=G'[V(C_5)]$, and
	\begin{equation}\label{equ-22}
		\chi(C_5)\geq\chi(C_4)-(2\ell-2)>0.
	\end{equation}
	
	We may assume that $T_2=p_1p_2\cdots p_{2\ell-2}\cdots p_m$, where $p_1=b_{t-1}$ and $p_2=b_t'$. By the cleanliness of $(V(C_5),T_2)$, $m\geq2\ell$ and $d_{G'}(p_i,V(C_5))\geq3$ for $1\leq i\leq2\ell-2$. Moreover, we have that
	\begin{equation}\label{equ-23}
		\mbox{$d_{G'}(b_{t-1},V(C_5))\geq4$ and $d_{G'}(a_t,V(C_5))\geq3$.}
	\end{equation}
	
	Let $p_{m+1}\in V(C_5)$ such that $p_m\sim p_{m+1}$, and let $c_{t-1}$ be a parent of $p_{m+1}$ in $M_{t-1}$ and $c_{t-2}$ be a parent of $c_{t-1}$ in $M_{t-2}$. Consequently, we have that $c_{i}\ne a_{i}$ for $i\in\{t-1,t-2\}$ and $a_{t-2}\not\sim c_{t-1}$ as $p_{m+1}\notin X$. By (\ref{equ-23}), $c_{i}\ne b_{i}$ for $i\in\{t-1,t-2\}$, $b_{t-1}$ is anticomplete to $\{c_{t-1},c_{t-2}\}$, $b_{t-2}\not\sim c_{t-1}$ and $a_t\not\sim c_{t-1}$. Since $(M_0',\cdots,M_t')$ is a week stable levelling of $G'$, and thus $\{a_{t-2},b_{t-2},c_{t-2}\}$ is a stable set.
	
	Choose the minimum $1\leq h\leq m+1$ such that $c_{t-1}\sim p_h$. Since $d_{G'}(c_j,V(C_5))\geq3$ for $1\leq j\leq2\ell-2$, $h\geq2\ell-1$. Let $P=c_{t-1}p_hp_{h-1}\cdots p_2$ (Notice that $p_2=b_t'$). Since $G'$ has no long odd hole, with the same argument in the proof of Lemma~\ref{L-7}, we have that $P$ is an induced path of odd length at least $2\ell-1$.
	
	If $c_{t-1}\sim a_{t-1}$, then $b_{t-1}\not\sim a_{t-1}$ by (\ref{equ-23}). Since $P$ is an induced path which has odd length at least $2\ell-1$, we have that the union of $P$, $p_2p_1b_{t-2}$, $c_{t-1}a_{t-1}a_{t-2}$ and an even ceiling-$a_{t-2}b_{t-2}$-path (in $\cup_{i=0}^{t}M_i$) is a long odd hole. Therefore, $c_{t-1}\not\sim a_{t-1}$.
	
	Suppose that $b_{t-1}\sim a_{t-1}$. If $c_{t-2}\sim a_{t-1}$, then the union of $P$ and $p_2p_1a_{t-1}c_{t-2}c_{t-1}$ is a long odd hole. If $c_{t_2}\not\sim a_{t-1}$, then the union of $P$, $p_2p_1a_{t-1}a_{t-2}$, $c_{t-1}c_{t-2}$ and an even ceiling-$a_{t-2}c_{t-2}$-path (in $\cup_{i=0}^{t}M_i$) is a long odd hole.
	
	Therefore, $\{a_{t-1},b_{t-1},c_{t-1}\}$ is a stable set, and we are done with the same argument in the proof of Lemma~\ref{L-7}.
	
	This completes the proof of Theorem~\ref{main-2}. \qed

\end{document}